\documentclass{amsart}
\usepackage{mathtools}

\numberwithin{equation}{section}
\newtheorem{thmA}{Theorem}

\newcommand{\R}{\mathbb{R}}
\newcommand{\N}{\mathbb{N}}

\newcommand{\T}{\mathbb{T}}

\newcommand{\B}{\mathbb{B}}

\newcommand{\pt}{\T[\![t]\!]}

  \usepackage[english]{babel}
     
  \usepackage{comment}
  \usepackage{all2021}                    %
  \usepackage{tree-commands4}             %
  \usepackage{bm}                         %
  
  \usepackage{overpic}                    %
  \usepackage{tikz-cd}                    %

  \usepackage[colorlinks, citecolor={blue}, linkcolor={{black}}, urlcolor={blue}]{hyperref}          %

  \calclayout
  
  \usepackage{enumerate}

\theoremstyle{definition}

\newcommand{\SolSet}{\mathcal{S}}
\usepackage{stmaryrd}
\newcommand{\pseries}[1]{\llbracket #1 \rrbracket}

\title{About finite differential tropical basis for linear ODE's}

\author{Stefano Mereta}
\email{stefano.mereta@cunef.edu}
\address{Departamento de Matemáticas, CUNEF Universidad, C. de Almansa 101, 28040 Madrid, Spain}
\author{Alejandro Vargas}
\email{alejandro@vargas.page}
\address{Mathematics Institute,
Zeeman Building,
University of Warwick,
Coventry CV4 7AL, UK}
\date{\today}


\begin{document} 

\begin{abstract}
We formulate several open questions regarding the tropicalization of linear ODEs, 
aiming primarily to develop methods for calculating the radius of convergence of their classical solutions. 
To this aim it is of foremost importance to characterize the classes of equations that admit a finite differential tropical basis, as introduced in \cite{ft22}.
Our initial exploration examines the second- and third-order cases.
\end{abstract}

		\maketitle  
 
\section{Introduction}
    \label{sec:Introduction}
Since at least the 1930s, algebraic techniques have been developed to solve (partial) differential equations expressed as polynomials in an indeterminate function and its derivatives. 
These efforts were carried out preeminently by J.~F.~Ritt and his student E.~R.~Kolchin (see \cite{ritt,kolchin}) and established \emph{differential algebra} as a new and effective toolkit to apply commutative algebra to study (partial) differential equations and their solution sets.

More recently, \emph{tropical geometry} (and algebra) emerged as an area of mathematics whose tools transform hard algebro-geometric problems into combinatorial ones.
These are more amenable for computational approaches.

In the last ten years, since \cite{grigoriev}, the use of tropical techniques for differential equations has been explored and its scope extended to differential equations with power series coefficients over arbitrary valued field (see \cite{framework}). 
While \cite{fundtriv, fundpartial, mg24} provide a proof of a fundamental theorem analogous to that of tropical geometry (see \cite[Theorem 3.2.5]{maclagansturmfels}) in full generality, 
the computational effectiveness of these techniques remains largely unexplored.

The finiteness of tropical bases (see \cite[Theorem 2.6.5]{maclagansturmfels}) is indispensable for implementing effective tropical algorithms in the classic algebro-geometric setting. 
When dealing with ideals generated by linear forms a tropical basis can always be chosen to consist of linear forms: the circuits of the valuated matroid associated to the linear space. 
Recently, \cite{ft22} introduced a notion of \emph{differential tropical basis} (DTB), where finiteness is defined up to the action of the derivative. 
However, the authors also demonstrate a class  of differential ideals generated by linear differential equations with constant coefficients for which there is no finite DTB of linear forms with respect to the trivial valuation.

This work explores further the finiteness of DTBs, 
taking the first steps towards characterizing this property for linear differential equations in the more general framework of non-trivial valuation. Our interest in DTB finiteness stems from its potential applications in computing the the radii of convergence of differential equations over valued fields via tropical methods (see \cite[Section 6]{mg24}).

We study in detail the second and third order case under the hypothesis that the roots of the characteristic polynomial have different non-zero valuation. 
We prove the following:

\begin{thmA}[See Theorem~\ref{thm:TropBasisOrder2}]
 Let $p$ be a prime and $F$ be the differential polynomial
\[ 
 F = (x-\alpha_1)(x- \alpha_2) = x^2 - (\alpha_1 + \alpha_2)x + \alpha_1 \alpha_2 \in \CC_p[\![t]\!]\{x\}.
\]
If $\nu_p(\alpha_1) \ne \nu_p(\alpha_2)$, then the set $\aset F$ is a DTB for the ideal~$[F]$. 
\end{thmA}

\begin{thmA}[See Theorem~\ref{thm:Determined}]
   Let $p$ be a prime and $ F = (x-\alpha_1)(x- \alpha_2)(x - \alpha_3)  \in \CC_p[\![t]\!]\{x\} $ a differential polynomial.
Set $u_i = \nu_p(\alpha_i)$ and $G^{(2)}_0 = \min(u_2 + u_3, u_2 + b_{1}, b_{2})$.  
Assume that $u_1  < u_2  < u_3 $. 
If the initial condition $(b_1,b_2)$ is not in $V^{\text{trop}}( G^{(2)}_0 )$, 
then the unique extension to a solution in $\bigcap_{n \in \NN}\SolSet(\trop(d^nF))$ is
\begin{align*} 
  b_n = (n-2) u_1 + G^{(2)}_0(b_1,b_2) - \nu_p(n!).
\end{align*}
\end{thmA}

The same methods used to prove that the set $\aset F$ is a DTB for the ideal~$[F]$ for the order 2 case also work for order 3, although the notation becomes quite cumbersome.
A generalization to order $r$ is straightforward, 
and we leave open finding theoretical machinery that proves this without lengthy calculations.
We formulate several other open questions in this note.

\subsection{Rudiments of differential and tropical algebra}
For comprehensive sources about differential algebra, we refer the reader to \cite{ritt, kolchin}. 
A \emph{differential ring} is a ring $R$ equipped with a derivation, i.e.\ an additive map $d_R : R \rightarrow R$ such that 
\[
    d_R(rs) = rd_R(s) + d_R(r)s
\]
holds for all $r,s \in R$. All along this work, the only differential ring we will consider is the ring of power series $K[\![t]\!]$ over a field $K$, equipped with the usual derivation $d/dt$. The \emph{differential algebra of differential polynomials} in one variable over $K[\![t]\!]$ is the $K[\![t]\!]$-algebra 
\[
    K[\![t]\!][x^{(j)}\mid j \in \N]
\]
equipped with the only derivation $d$ extending $d/dt$ and such that $dx^{(j)} = x^{(j+1)}$. 
We will denote this differential algebra as $K[\![t]\!]\{x\}$. 
An ideal $I \subseteq K[\![t]\!]\{x\}$ is \emph{differential} if $df \in I$ for all $f \in I$. 
Given a differential polynomial $f$, the \emph{differential ideal generated by} $f$ is the smallest differential ideal containing $f$, 
i.e.~the (algebraic) ideal generated by the family $\{d^n f \mid n \in \N\}$. 
We will denote this ideal as $[f]$.

For fundamentals of valuations and tropicalization, 
we refer the reader to \cite[Section~1.7 and~1.8]{JeffNotes} in this volume. 
We just recall the definition of the semirings we will use in the following, and that of valuation.
A \emph{semiring} $(S, \oplus, \odot)$, is an algebraic structure satisfying the same axioms as a ring, without the requirement that $(S, \oplus)$ is a group. If $(S \setminus \{0_S\}, \odot)$ is a group, we will say that $S$ is a \emph{semifield}. A semiring is idempotent if for every element $a \in S$ it holds that $a \oplus a = a$. All the semirings appearing in this work are idempotent, the following are the relevant examples:
\begin{itemize}
    \item The semifield of (rank 1) tropical numbers $(\T, \oplus, \odot) :=(\R \cup \{\infty\}, \min, +)$;
    \item The semifield $\B$ of boolean numbers i.e.\ the set $\{0,\infty\}$ with the semiring structure inherited as a subset of $\T$;
    \item The semifield of rank $n$ tropical numbers $(\T^{(n)}, \oplus, \odot) :=(\R^n \cup \{\infty\}, \min_{\text{lex}}, +)$ where $\min_{\text{lex}}$ is the minimum with respect to the lexicographic order on $\R^n$.
\end{itemize}
\begin{de}
A \emph{valuation} on a field $K$ is a map $\nu : K \to \T^{(n)}$ that satisfies the following three axioms, for all $a,b\in K$:
\begin{enumerate}
    \item\label{axiom: valuatioun1} $\nu^{-1}(\infty)=\{0\}$;
    \item\label{axiom: valuatioun2} $\nu(ab) = \nu(a)\odot \nu(b)$;
    \item\label{axiom: valuatioun3} $\nu(a+b)\ge \nu(a)\oplus\nu(b)$.
\end{enumerate}
\end{de}
Now, we briefly present the notions of tropical differential algebra required in the following, as introduced in \cite{framework}. For an element $A \in K[\![t]\!]$ we denote as $\textup{Supp}(A)$ its support. If the field $K$ is endowed with a valuation $v_K : K \rightarrow \T$, then the assignment $v :K[\![t]\!] \rightarrow \T^{(2)} $ defined as 
\[
   A = \sum_{i=0}^\infty a_it^i \mapsto (\min \textup{Supp}(A), v_K(a_{\min \textup{Supp}(A)}))
\]
is a valuation. If $v_K$ is non-trivial $v$ is a rank 2 valuation, otherwise it is of rank 1. 
The tropicalization of differential polynomials takes place via the valuation $v$: given a differential polynomial $f \in K[\![t]\!]\{x\}$ its tropicalization $\trop_v(f)$ is the element of $\T^{(2)}[x^{(j)}\mid j \in \N]$ obtained by applying $v$ to the coefficients of $f$. For a differential ideal $I \subset K[\![t]\!]\{x\}$ its tropicalization $\trop_{v}(I)$ is
\[
    \trop_{v}(I) = \{\trop_{v}(f) \mid f \in I\}.
\]

Analogous to the ring-theoretic construction, tropical power series form a semiring $\TT \pseries t$.
This becomes a differential semiring (for details about them see \cite[Section 2.4]{framework}) by endowing it with the differential $(\frac{d}{dt})_{v_K}$ defined by the rule:
\[
    \left ( \frac{d}{dt} \right )_{v_K}(t^n) = 
    \begin{cases}
        v_K(n)t^{n-1} & \text{if $n \ge 1$} \\        \infty & \text{if $n = 0$}.
    \end{cases}
\]
We denote the differential semiring $(\T[\![t]\!], \left ( \frac{d}{dt} \right )_{v_K}) $ as $\T[\![t]\!]_{v_K}$.

In contrast with classic tropical geometry,
power series in $K[\![t]\!]$ are not tropicalized with respect to the same map used for equations: let $A \in K[\![t]\!]$, then its tropicalization $\trop(A)$ is obtained by applying $v_K$ to its coefficients, thus obtaining an element of $\T[\![t]\!]_{v_K}$. We denote this assignment $\trop_{v_K}: K[\![t]\!] \rightarrow \T[\![t]\!]_{v_K}$ and drop the subscript if the situation makes it clear. It follows from the definitions that 
\[
    \left ( \frac{d}{dt} \right )_{v_K}(\trop_{v_K}(A)) = \trop_{v_K}(dA)
\]
for every $A \in K[\![t]\!]$.

Finally, in order to evaluate an element $f$ of $\T^{(2)}[x^{(j)}\mid j \in \N]$ at a tropical power series $B \in \T[\![t]\!]$, we introduce the homomorphism of semiring $w:  \T[\![t]\!]_{v_K} \rightarrow \T^{(2)}$ defined in a similar fashion to the valuation $v$: 
\[
     B = \sum_{i=0}^\infty b_it^i \mapsto (\min \textup{Supp}(B), b_{\min \textup{Supp}(B)}).
\]

The evaluation of $f \in \T^{(2)}[x^{(j)}\mid j \in \N]$ at a tropical power series $B$ is the element of $\T^{(2)}$ obtained by plugging $(d/dt)^j_{v_K}(B)$ for $x^{(j)}$. 
We say that $B$ is a \emph{tropical solution} to $f$ if the resulting minimum in $\T^{(2)}$ is attained at least twice.

\subsection{Linear homogeneous ODEs with constant coefficients}
    \label{sub:LinearhomogeneousODE}

Let $K$ be an algebraically closed field. Recall that a linear homogeneous ODE of order $r$ with constant coefficients $c_i \in K$ and unknown function $x(t)$ is expressed as:
\begin{align}
    \label{eq:homogeq}
    F = x^{(r)} + c_{r-1} x^{(r-1)} + \dots + c_1x^{(1)} + c_0x = 0,
\end{align}
where the $c_i$ are in $K$ and $x^{(i)}$ is the derivative $\frac{d^i}{dt^i} x(t)$.
The solution set $\SolSet(F)$ is a linear space of dimension~$r$,
determined by the associated characteristic polynomial:
\[
    \chi(z) = z^{r} + c_{r-1} z^{r-1} + \dots + c_1z + c_0.
\]
Let $\kappa$ be the number of roots of $\chi(z)$ and let $\aset{\alpha_s}_{s \in [\kappa]}$ be the set of its roots. Denote as $m_s$ the multiplicity of $\alpha_s$. 
The functions $\{ t^{j-1}\exp(\alpha_s t) \mid j = 1, \dots , m_s\}_{s \in [\kappa]}$, defined by the power series expansion 
\[ t^{j-1}\exp(\alpha_s t) = \sum_{i = 0}^\infty \frac{\alpha_s^i}{i!}t^{i+j-1} \]
are solutions to Equation~\eqref{eq:homogeq} and the set $\{ t^{j-1}\exp(\alpha_s t) \mid j = 1, \dots , m_s\}_{s \in [\kappa]}$ is a basis for $\SolSet(F)$.
We assume that $c_0 \ne 0$, so that $0$ is not a root of $\chi(z)$, and that all the roots are pairwise distinct, i.e. $\kappa = r$ and $m_s = 1$ for all $s \in [r]$: under these hypothesis $\aset{ \exp(\alpha_s t)  }_{s \in [r]}$ is a basis for $\SolSet(F)$.

\subsection{Tropicalization of ODEs under trivial valuation}
    \label{sub:TropicalizationOfODEs}
Under the trivial valuation, the tropicalization of $\SolSet(F)$ yields a matroid $M$ over $\NN$~\cite{abf+23a}
(see~\cite{bdk+13} for the relevant axioms on infinite matroids). 
To compute $M$, 
we construct an $r \times \NN$ matrix using the power series coefficients of our basis functions~$\exp(\alpha_s t)$.
The bases of $M$ are sets $I \subset \NN$ indexing non-zero maximal minors.
Because column rescaling does not alter the underlying matroid, 
we factor out the $1/i!$ terms and analyze the following matrix:
    \[A = 
       \begin{pmatrix}
1 & \alpha_1 & \alpha_1^2 & \alpha_1^3 & \alpha_1^4 & \dots \\
1 & \alpha_2 & \alpha_2^2 & \alpha_2^3 & \alpha_2^4 & \dots \\
 \vdots&         & \vdots      &          &       \vdots        & \\
\\
1 & \alpha_r & \alpha_r^2 & \alpha_r^3 & \alpha_r^4 & \dots \\
\end{pmatrix} 
    \]
    If all maximal minors are non-zero, $M$ is the finitary rank $r$ uniform matroid $U_{r, \NN}$. 
    Over fields of characteristic zero, this generic behavior follows from the properties of generalized Vandermonde determinants. 

    In the non-generic case, suppose $I \in {\mathbb N \choose r}$ indexes a vanishing minor. 
    Let $k_1 < k_2 < \dots < k_r$ denote the elements of $I$. 
    Factoring out the lowest power yields
\begin{align}
\begin{vmatrix}
\alpha_1^{k_1} & \alpha_1^{k_2} & \cdots & \alpha_1^{k_r}\\
\vdots & \vdots & \ddots &  \vdots  \\
\alpha_r^{k_1} & \alpha_r^{k_2} & \cdots & \alpha_r^{k_r} \\
\end{vmatrix} =
\alpha_1^{k_1} \alpha_2^{k_1}  \cdots \alpha_r^{k_1} \begin{vmatrix}
1& \alpha_1^{k_2-k_1} & \cdots & \alpha_1^{k_r - k_1}\\
\vdots & \vdots & \ddots &  \vdots  \\
1 & \alpha_r^{k_2-k_1} & \cdots & \alpha_r^{k_r-k_1} \\
\end{vmatrix}, 
\end{align}
demonstrating that shifting $I$ preserves the vanishing of the corresponding minors (as $\alpha_i \neq 0$ for all $i$):

\begin{prop} 
    \label{prop:ShiftOperatorAndMatroids}
     Let $\sigma : \ZZ \to \ZZ$ denote the shift map $i \mapsto i+1$.
     For any integer $q$ such that $\sigma^q(I) \subset \NN$, the shifted index set $\sigma^q(I)$ yields a vanishing minor if and only if $I$ does. 
\end{prop}

This observation motivates the following question:

\begin{que} 
    \label{que:GroupAction}
 Can the infinite matroid $M$ on the ground set $\NN$ associated with a linear homogeneous ODE of order~$r$, 
 be described using only a finite matroid and a group action?  
\end{que}

\subsection{Differential tropical basis and radius of convergence}
    \label{sub:NontrivialValuationRadiusOfConvergence}

    Following \cite[Section~6]{mg24}, 
    we define the \emph{radius of convergence with respect to $c \in \RR_{>1}$} 
    of a tropical power series $\bigoplus_{i = 0}^\infty b_i t^i \in \TT \pseries t$ as
        \[
             \rho_c \left( \sum_{i = 0}^\infty b_i t^i \right) 
             := \sup \aset{r \in [0, \infty) \mid 
             \lim_{i \to \infty} c^{-b_i} r^i = 0}.  
        \]
    Here the convention is that $c^{-b_i} = 0$ when $b_i = \infty$.

    In contrast, let $K$ be an uncountable algebraically closed field of characteristic zero equipped with a nontrivial valuation $v_K : K \to \TT$ admitting a section.
    We choose $c \in \RR_{>1}$ and define the norm $|x|_{K,c} = c^{-v_K(x)}$.
    The \emph{radius of convergence} of a power series $\sum_{i = 0}^\infty a_i t^i \in K \pseries t$ is given by 
        \[
             \rho_c \left( \sum_{i = 0}^\infty a_i t^i \right) 
             := \sup \aset{r \in [0, \infty) \mid 
             \lim_{i \to \infty} |a_i|_{K,c} r^i = 0}.  
        \] 
    A direct consequence of the main theorem in~\cite{mg24} is that computing the convergence radius commutes with tropicalization.

    \begin{prop} 
        For any $F \in K\pseries t \{x \}$, we have:
        \begin{align} 
            \label{eq:ConvergenceRadius}
        \aset{\rho_c(A) \suchthat A \in \SolSet(F) }  = \aset{\rho_c(B) \suchthat B \in \SolSet( \trop([F]) )}.
        \end{align} 
    \end{prop}

    The right-hand side of Equation~\eqref{eq:ConvergenceRadius} is tropical, 
    and thus amenable to computation via polyhedral methods and $(\min, +)$-algebra algorithms.
    In particular, \cite[Example~6.4]{mg24} demonstrates that having a finite DTB at hand facilitates this computation:
    \begin{de}
        Let $I \subseteq K[\![t]\!]\{x\}$ be a differential ideal. As introduced in  \cite{ft22}, a \emph{differential tropical basis} (DTB) for $I$ is a set of differential polynomials $\{F_\lambda\}_{\lambda \in \Lambda} \subset I$ for some index set $\Lambda$ such that 
\[
    \SolSet(\trop(I)) = \bigcap_{\lambda \in \Lambda} \bigcap_{n \in \NN} \SolSet(\trop(d^nF_\lambda)).
\]
\end{de}

This notion is the differential analogue to the classic notion of tropical basis for an ideal (see \cite[Definition 2.6.4]{maclagansturmfels} for the definition of tropical basis).

The discussion above motivates the following question:
   
\begin{que} 
    \label{que:FiniteBasis}
    Which linear differential polynomials $F \in K\pseries t \{x \}$ admit a \emph{finite} differential tropical basis for $[F]$?
\end{que}

     The matter is subtle: 
     \cite[Section~4.2]{ft22} shows that, generically, 
     a second-order linear ODE does not have a finite DTB of linear forms 
     (this does not rule out the existence of a finite tropical basis of higher degree, \emph{a priori}). 
    On the other hand, 
    in Subsection~\ref{sub:padicvaluation} we identify two sufficient conditions on the equation for a finite DTB of linear forms to exist.
    In the spirit of Question~\ref{que:GroupAction}, one can ask:

\begin{que} 
    \label{que:AdditionalAction}
     Can we define a notion of differential tropical basis that incorporates an additional group action, beyond the differential one, so that all linear ODEs have a finite differential tropical basis of linear forms? 
\end{que}

\section{Second order}
    \label{sub:OrderTwoExample}

\subsection{Trivial valuation}
    \label{sub:OrderTwoTrivialValuation}
    Consider Equation~\eqref{eq:homogeq} for the case~$r=2$, on the complex numbers with a trivial valuation. 
    Let $\alpha$ and $\beta$ denote the roots of the characteristic polynomial $\chi(z)$. 
    Assuming $\alpha \ne \beta$, the matroid $M$ is represented by the matrix
\[
A = 
       \begin{pmatrix}
1 & \alpha & \alpha^2 & \alpha^3 & \alpha^4 & \dots \\
1 & \beta & \beta^2 & \beta^3 & \beta^4 & \dots
\end{pmatrix} 
\]
If all maximal minors are non-zero, $M$ is the infinite uniform matroid $U_{2, \NN}$. 
Otherwise, we choose indices $k<l$ such that the minor formed by the $k$-th and $l$-th columns vanishes, and the difference $l - k$ is minimal. 
Factoring out the lowest powers yields:
\[
\begin{vmatrix}
\alpha^k & \alpha^l\\
\beta^k & \beta^l
\end{vmatrix} =
\alpha^k \beta^k \begin{vmatrix}
1 & \alpha^{l-k}\\
1 & \beta^{l-k}
\end{vmatrix}
= \alpha^k \beta^k (\beta^{l-k} - \alpha^{l-k}) =0. 
\]
Recall that $\alpha, \beta$ are non-zero, thus $\beta^{l-k} - \alpha^{l-k} = 0$ means that $\omega = \beta/\alpha$ is a $(l-k)$-th root of unity. 
So $\beta = \omega \alpha$ and we rescale the columns of $A$ to obtain the following matrix 
    \[\tilde A = 
       \begin{pmatrix}
1 & 1 & 1 & 1 & 1 & \dots & 1 & 1 & 1 & \dots \\
1 & \omega & \omega^2 & \omega^3 & \omega^4 & \dots & \omega^{l-k-1} & 1 & \omega & \dots 
\end{pmatrix} 
    \]

Consequently, a set $\{a,b\}$ is a base if and only if $a \not \equiv b \mod (l-k)$.
The base polytope $P_M$ associated with $M$ possesses infinitely many vertices  $\mathbf  e_{\aset{a,b}}$;
however, $P_M$ is determined by the polytope associated with the finite uniform matroid $U_{2, l - k}$
plus an action on $\RR^\NN$ by the permutation group $H \subset S_\NN$, generated by all the transpositions of the form $(i j)$ satisfying $i \equiv j \mod (l-k)$. 
The number of orbits is finite, having canonical representatives 
$\aset{\bar a, \bar b}$ where $\bar a, \bar b \in \aset{0, \dots, l - k - 1}$.

Regarding the flats: proper flats have rank 1 and there are $(l-k)$-many of them, namely the $l-k$ equivalence classes of $\mathbb N$ modulo $l-k$. 
Thus, remarkably, the Bergman fan of $M$  consists of finitely many cones, despite the infinite ground set. 
See Figure~\ref{fig:FiniteRays} below.

\begin{figure}[htpb]
  \label{fig:FiniteRays}
    \centering
    \begin{tikzpicture}[scale=0.8, >=stealth]
        \filldraw (0,0) circle (1.5pt) node[anchor=north east] {$\mathbf{0}$};

        \draw[->, thick] (0,0) -- (0, 2) node[above] {$\mathbf{e}_{[0]}$};
        
        \draw[->, thick] (0,0) -- (1.9, 0.6) node[right] {$\mathbf{e}_{[1]}$};
        
        \draw[->, thick] (0,0) -- (1.1, -1.6) node[below right] {$\mathbf{e}_{[2]}$};

        \draw[dotted, thick] (0,-1) arc[start angle=-90, end angle=-150, radius=1];

        \draw[->, thick] (0,0) -- (-1.9, 0.6) node[left] {$\mathbf{e}_{[l-k-1]}$};

    \end{tikzpicture}
    \caption{The Bergman fan of $M$. Despite the infinite ground set, the fan consists of exactly $l-k$ rays corresponding to the equivalence classes of $\mathbb{N} \pmod{l-k}$.}
    \label{fig:BergmanFan}
\end{figure}
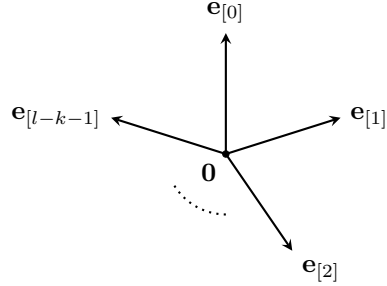

It is natural at this point to ask about the Boolean dimension of the tropicalization of the solution set: this is the cardinality of a minimal system of generators.
Interpreting the members of the tropical solution set as indicator vectors of subsets of $\NN$, 
the Boolean dimension counts the number of inclusion minimal sets, 
thus it can be infinite.
In this example, we have
\[\dim_\BB(\Trop( \rowspan S)) = l - k,\]
as the support of the generators is disjoint.
Note that classically the dimension is 2, thus showing that dimension can increase arbitrarily under tropicalization.

Last, if $\alpha = \beta$, then the solutions are of the form
  \begin{align*}
        \exp(\alpha t), \quad t\exp(\alpha t),
    \end{align*}
and so we get the matrix:
 \[S = 
       \begin{pmatrix}
1 & \alpha & \frac 1 2 \alpha^2 & \frac 1 3 \alpha^3 & \frac 1 4 \alpha^4 & \dots \\
0 & 1      & \alpha & \alpha^2 & \alpha^3  & \dots \\
\end{pmatrix}. 
    \]
    It is straightforward to see that this corresponds to the infinite uniform matroid $U_{2,\NN}$.

\subsection{$p$-adic valuation}
    \label{sub:padicvaluation}
We revisit the example from~\cite[Section~4.2]{ft22}, which demonstrates that a generic second-order linear equation with constant coefficients lacks a finite DTB of linear forms in the trivially valued case. 
Shifting to the non-trivially valued case 
we show that if the roots of the characteristic polynomial have distinct valuations, 
a finite DTB is guaranteed to exist.

Let $(\CC_p, \nu_p)$ denote the valued field of $p$-adic complex numbers for a prime number $p$, and consider the differential polynomial
\[ 
 F = (x-\alpha_1)(x- \alpha_2) = x^2 - (\alpha_1 + \alpha_2)x + \alpha_1 \alpha_2 \in \CC_p[\![t]\!]\{x\}.
\]
As recalled above, the classical solutions are $\exp(\alpha_1t)$ and $\exp(\alpha_2t)$. Let $u_1 = \nu_p(\alpha_1)$ and $u_2 = \nu_p(\alpha_2)$ denote the valuations of the roots and set $c=p$ as its customary for the $p$-adic norm. It follows from the fact that the radius of convergence of the $p$-adic exponential is $p^{-\frac{1}{p-1}}$ that the radius of convergence of these two series is:
\begin{align*}
    \rho(\exp(\alpha_1t)) &= p^{-\frac{1}{p-1} + u_1} \\
    \rho(\exp(\alpha_2t)) &= p^{-\frac{1}{p-1} + u_2}
\end{align*}

Assuming $u_1 < u_2$, 
the properties of valuations imply that $\alpha_1 \ne 0$ and  $\nu_p(\alpha_1 + \alpha_2) = u_1$, with the case $\alpha_2 = 0$ corresponding to~$u_2 = \infty$.

We aim to characterize the tropical power series $B = \sum_{i \in \NN} b_i t^i  \in \pt_{\nu_p}$ solving the system $\{\trop (d^n F) \mid n \in \N\}$. 
Since $F$ has order $r=2$, the initial coefficients $b_0$ and $b_1$ completely determine the series. 
By rescaling the solution, we may assume without loss of generality that $b_0 = 0$,
and treat $b_1$ as a free parameter to compute the subsequent terms. 
Evaluating $\trop (d^nF)$ at $B$ yields the following tropical vanishing condition for all $n \in \N$:
\begin{align}
    \label{eq:Degree2Example}
    (0,b_n + u_1 + u_2 + \nu_p(n!)) \oplus 
    (0,b_{n+1} + u_1 +\nu_p((n+1)!)) \oplus 
    (0,b_{n+2} + \nu_p((n+2)!)).
\end{align}
If $u_2 = \infty$, then Equation~\eqref{eq:Degree2Example} reduces to the vanishing of
\begin{align}
    (0,b_{n+1} + u_1 +\nu_p((n+1)!)) \oplus 
    (0,b_{n+2} + \nu_p((n+2)!)).
\end{align}
This immediately yields the recurrence relation $b_{n+2} = b_{n+1} + u_1 - \nu_p(n+2)$ for all $n \ge 1$.

It remains to consider the case where $u_2 < \infty$.
Setting $n=0$ in Equation~\eqref{eq:Degree2Example} requires the following expression to tropically vanish: \[(0, u_1 + u_2 ) \oplus 
    (0,b_{1} + u_1 ) \oplus 
  (0,b_{2} + \nu_p(2!)).\] 
This is equivalent to defining the tropical polynomial 
\[  
    G_0 = \min(u_1 + u_2,   u_1 +  b_{1}  , b_2 + \nu_p(2!) ) \in \TT[b_1, b_2],
\]
and requiring the point $(b_1,b_2)$ to lie in the tropical variety $V^{\trop}(G_0) \subset \TT^2$. 
In other words, the minimum in $G_0$ must be achieved at least twice. 
This requirement yields three distinct cases,  depending on the value of the parameter $b_1$.
In the following boxed equations we follow the first case as it branches along.
\begin{enumerate}[(1)]

    \item \label{item:b1LessThanU2} if $b_1 < u_2$: the minimum value is $b_1 + u_1$, achieved by the second and third terms of $G_0$; forcing $\boxed{b_2 = b_1 + u_1 - \nu_p(2!)}$.
    \item \label{item:b1GreaterThanU2} if $b_1 > u_2$: the minimum is $u_1 + u_2$, achieved by the first and third terms; forcing $b_2 = u_1 + u_2 - \nu_p(2!)$.
    \item \label{item:b1EqualToU2} if $b_1 = u_2$: the first two terms coincide: $u_1 + u_2 = b_1 + u_1$. 
      So $G_0$ vanishes tropically if and only if $b_2 \ge u_1 + u_2 - \nu_p(2!)$.
\end{enumerate}

Next, analogously as above, setting $n=1$ in Equation~\eqref{eq:Degree2Example} yields:
\[  G_1 = \min
    (b_1 + u_1 + u_2,
    b_{2} + u_1 + v_p(2!), 
    b_{3} + \nu_p(3!)) \in \TT[b_2, b_3].\]
To compute the value of $b_3$ that makes $G_1$ tropically vanish, 
the cases are as follows: 
\begin{enumerate}[(1)]
    \item if $b_1 < u_2$: recall that $b_2 = b_1 + u_1 - v_p(2!)$, 
      so $ G_1 = \min (b_1 + u_1 + u_2, b_1 + 2  u_1, b_{3} + \nu_p(3!))$.
        As $u_1 < u_2$, the minimum is $b_1 + 2u_1$, achieved by the second and third terms of $G_1$; forcing $\boxed{b_3 = b_1 + 2u_1 - \nu_p(3!)}$.
    \item if $b_1 > u_2$: recall that $b_2 = u_1 + u_2 - v_p(2!)$, so $ G_1 = \min (b_1 + u_1 + u_2, u_2 + 2  u_1 , b_{3} + \nu_p(3!))$.
    As $u_1 < u_2 < b_1$, the minimum is $u_2 + 2u_1$, achieved by the second and third terms of $G_1$; forcing $b_3 = u_2 + 2u_1 - \nu_p(3!)$.

    \item $b_1 = u_2$: recall that $u_1 + u_2 = b_1 + u_1$, so $G_1 = \min (u_1 + 2u_2, b_2 +  u_1 +v_p(2!), b_{3} + \nu_p(3!))$.
      Also recall that $b_2 \ge u_1 + u_2 - \nu_p(2!)$, so the second term is not smaller than the first term.
        Now we are in an identical position as the $n=0$ case, 
        yielding three sub-cases:
        \begin{enumerate}[(3.1)]
            \item  $b_2 < 2u_2 - v_p(2!)$: proceeds analogously to case~\eqref{item:b1LessThanU2}.
              From $b_2 < 2u_2 - v_p(2!)$ we get $b_2 + v_p(2!)  + u_1 < 2u_2 + u_1$
              forcing $b_3 = b_2 + \nu_p(2!) + u_1  - \nu_p(3!)$.
            \item $b_2 > 2u_2 -v_p(2!)$: proceeds analogously to case~\eqref{item:b1GreaterThanU2}.
              From $b_2 > 2u_2 - v_p(2!)$ we get $b_2 + u_1 + v_p(2!) > 2u_2 + u_1$
              forcing $b_3 = 2u_2 + u_1 - v_p(3!)$.
            \item $b_2 = 2u_2 - v_p(2!)$: We obtain $G_1 = \min (u_1 + 2u_2, u_1 + 2u_2, b_{3} + \nu_p(3!))$. 
        The only constraint on $b_3$ is the inequality $b_3 \ge u_1 + 2u_2 - \nu_p(3!)$.
        \end{enumerate}
\end{enumerate}
By induction, one can prove that in cases 1 and 2:
\begin{enumerate}[(i)]
    \item Case $b_1 < u_2$: 
      assume $b_k = b_1 + (k-1)u_1 - \nu_p(k!)$.  
        Evaluating, we have $ G_{k-1} = \min (b_1 + (k-1) u_1 + u_2, b_1 + k u_1 , b_{k+1} + \nu_p((k+1)!))$.
    As $u_1 < u_2$, the second term is the unique minimum among the first two; 
  thus $G_{k-1}$ tropically vanishes if and only if
    $\boxed{b_{k+1} = b_1 + ku_1 - \nu_p((k+1)!)}$.
    \item $b_1 > u_2$: similar reasoning leads to the closed-form $b_{k+1} = u_2 + ku_1 - \nu_p((k+1)!)$ in this case.
\end{enumerate}
Both of these families of solutions yield a tropical radius of convergence equal to $\rho(\exp(\alpha_1t))$. 

Case 3 exhibits unstable behaviour. 
On the  one hand, at each step  one can keep choosing $b_i = i u_2 - \nu_p(i!)$ to build a solution with a tropical radius of convergence equal to  $\rho(\exp(\alpha_2t))$.
On the other hand, at each step the parameter $b_i$ only needs to satisfy the constrain $b_i \ge u_1 + (i-1)u_2 - \nu_p(i!)$, 
thus a choice distinct than $b_i = i u_2 - \nu_p(i!)$ reverts the sequence to a solution with a radius of convergence $\rho(\exp(\alpha_1t))$

These branching solutions are fully determined by a branching index $n \in \NN \cup \{\infty\}$ and a parameter $b \in \TT$ such that $b \ge u_1 + (n-1)u_2-\nu_p(n!)$ when $n > 1$. 
Specifically, a sequence that eventually behaves as in Case 1, namely $b<nu_2-\nu_p(n!)$, is given by
\begin{align*} 
    b_i & = i u_2 - \nu_p(i!) \text{ for all } i < n, \\
    b_n & = b, \\
    b_{n+i} &= b + \nu_p(n!) + iu_1 - \nu_p((n+i)!) \text{ for all } n+i > n.
\end{align*}
A sequence that eventually behaves as in Case 2, namely $b>nu_2 - \nu_p(n!)$, is given by:
\begin{align*} 
    b_i & = i u_2 - \nu_p(i!)  \text{ for all } i < n, \\
    b_n & = b, \\
    b_{n+i} &= nu_2 + \nu_p(n!) + iu_1 - \nu_p((n+i)!)  \text{ for all } n+i > n.
\end{align*}
Moreover, all the solutions we have described are realizable, 
that is, they are the tropicalizations of classical solutions. 
Indeed, consider the following convenient elements:
\begin{align*} 
  \psi_n =& \frac{n! p^b}{\alpha_1^n} \exp(\alpha_1 t ) \\
  =& \frac{n!p^b}{\alpha_1^n} +  \frac{n!p^b}{\alpha_1^{n-1}} t + \frac 1 2 \frac{n!p^b}{\alpha_1^{n-2}} t^2 + \dots \\
  \varphi_n =& \frac{\exp(\alpha_1 t )}{\alpha_1^n} - \frac{\exp(\alpha_2 t )}{\alpha_2^n}  = \sum_{i =0}^\infty \frac 1 {i!} (\alpha_1^{i-n} - \alpha_2^{i-n}) t^i \\
  =& \left(\frac 1 {\alpha_1^n} - \frac 1 {\alpha_2^n} \right ) 
            +   \left(\frac 1 {\alpha_1^{n-1}} - \frac 1 {\alpha_2^{n-1}} \right ) t + \dots  
            + \frac 1 {(n-1)!} \left(\frac 1 {\alpha_1} - \frac 1 {\alpha_2} \right)  t^{n-1}  \\
            &+ \frac 1 {(n+1)!}(\alpha_1 - \alpha_2)t^{n+1} + \dots
\end{align*}
We write $[A]_i$ for the $i$-th term of the series $A$.
Recall $b + \nu_p(n!) \ge u_1 + (n-1)u_2$, thus:
\begin{align*} 
  [\trop (\psi_n)]_i &= b + \nu_p(n!) + (i-n)u_1 - \nu_p(i!) \\ 
                       &\ge iu_2 + (n-i-1)(u_2-u_1) - \nu_p(i!)     \\
  [  \trop( p^{nu_2} \varphi_n)]_i & =  
\begin{cases}
iu_2 - \nu_p(i!) & \text{if } i < n, \\
\infty & \text{if } i = n, \\
nu_2 + (i - n)u_1 - \nu_p(i!) & \text{if } i > n.
\end{cases}
\end{align*}
 
Note above that $(u_2 - u_1)$ is positive. 
Given parameters $n$ and $b$, we have two cases. 
If $b < nu_2 - \nu_p(n!)$, there is a generic element $\epsilon$ of valuation $0$ such that 
\begin{align*}
[\trop (\exp(\alpha_2 t) + \epsilon \psi_n)]_i &= \min( iu_2 - \nu_p(i!), b + \nu_p(n!) + (i-n)u_1 - \nu_p(i!)) \\
&= \begin{cases}
iu_2 - \nu_p(i!) & \text{if } i < n, \\
b & \text{if } i = n, \\
b + \nu_p(n!) + (i - n)u_1 - \nu_p(i!) & \text{if } i > n.
\end{cases}
\end{align*}
If  $b \ge  nu_2 - \nu_p(n!)$, set $b' = b - ( nu_2 - \nu_p(n!)) \ge 0$.
Then, 
again up to multiplying by a generic element $\epsilon$ of valuation $0$, 
we can realize a branching that eventually behaves as Case 2 by: 
\begin{align*} 
  [\trop (p^{b'}\exp(\alpha_2 t) + \epsilon p^{nu_2} \varphi_n)]_i =
 \begin{cases}
iu_2 - \nu_p(i!) & \text{if } i < n, \\
b & \text{if } i = n, \\
nu_2 + (i - n)u_1 - \nu_p(i!) & \text{if } i > n.
\end{cases} 
\end{align*}
Finally, $\exp(\alpha_2t)$ realizes the unstable solution.
Thus: 

\begin{thm} 
    \label{thm:TropBasisOrder2}
 Let $p$ be a prime and $F$ be the differential polynomial
\[ 
 F = (x-\alpha_1)(x- \alpha_2) = x^2 - (\alpha_1 + \alpha_2)x + \alpha_1 \alpha_2 \in \CC_p[\![t]\!]\{x\}.
\]
If $\nu_p(\alpha_1) \ne \nu_p(\alpha_2)$, then the set $\aset F$ is a DTB for the ideal~$[F]$. 
\end{thm}

On the other hand, the following question remains open:

\begin{que} 
    \label{que:SameValuations}
    What are the possible behaviours when the roots $\alpha_1, \alpha_2$ have the same valuation?
    Is there an ideal~$[F]$ with a finite DTB that consists of more than one element?
\end{que}

We also wish to mention the following question.

\begin{que}
   Under trivial valuation, 
   is having a tropical basis, 
   in the sense of Fink and Toghani, 
   equivalent to having an associated polyhedral complex with finitely many cones? 
\end{que}

By what we have seen in this section, the assumption of trivial valuation is needed to pose the question, 
otherwise roots with different valuation give rise to infinitely many polyhedral cells.

\section{Order 3}
    \label{sub:OrderThreeExample}

\subsection{Trivial valuation}
    \label{sub:OrderThreeTrivialValuation}
We now consider Equation~\eqref{eq:homogeq} for~$r=3$ and describe the possible matroids obtained from the matrix~$A$. 
On the one hand, we still have the two extremes seen in the order~2 case, 
namely a periodic matrix giving rise to a finite-up-to-symmetry matroid, versus the uniform matroid $U_{3, \NN}$.
On the other hand, we provide several examples of \emph{intermediate} behaviors, 
and relate the vanishing of maximal minors to the secants of the monomial curve $t \mapsto (t^n, t^m)$ passing through $(1,1)$.
We leave open the question of fully classifying all the matroids that arise in this case. 

Let $\alpha, \beta, \gamma$ denote the roots of the characteristic polynomial~$\chi(z)$. 
Assuming the roots are distinct, the matroid $M$ is represented by the matrix
    \[S = 
       \begin{pmatrix}
1 & \alpha & \alpha^2 & \alpha^3 & \alpha^4 & \dots \\
1 & \beta & \beta^2 & \beta^3 & \beta^4 & \dots\\
1 & \gamma & \gamma^2 & \gamma^3 & \gamma^4 & \dots
\end{pmatrix}.
    \]
Again, for generic choices of $\alpha$, $\beta$, and $\gamma$ we get the uniform matroid  $U_{3, \NN}$. 
Otherwise, we list some examples.

\subsubsection{There is a big rank-1 flat}
Suppose there is a rank-1 flat with cardinality greater than 1.
This means there are two columns that are linearly dependent, say for $m<n$:
    \begin{align*} 
 \begin{psmallmatrix}
\alpha^m & \alpha^n \\
\beta^m & \beta^n \\
\gamma^m & \gamma^n \\
\end{psmallmatrix} 
    \end{align*}
The vanishing of the maximal minors yields the equations $(\beta/\alpha)^{n-m} - 1 = 0$ and $(\gamma/\alpha)^{n-m} - 1 = 0$. 
Thus, the ratios are once again roots of unity, yielding a periodic matrix $A$.

\subsubsection{Small rank-1 flats, some big rank-2 flats}
    Suppose there is a rank-2 flat with cardinality greater than 2.
    This means there is at least one vanishing minor, and we choose indices $k<l<m$ for such a vanishing minor. 
    Via row and column operations it can be shown that the vanishing of the minor of $S_{k,l,m}$ is equivalent to the vanishing of 
    \begin{align*} 
\label{eq:DeterminantOrder3}
 \begin{vsmallmatrix}
1& 1 & 1\\
0&(\beta/\alpha)^{l-k} - 1 & (\beta/\alpha)^{m-k} - 1\\
0&(\gamma/\alpha)^{l-k} - 1& (\gamma/\alpha)^{m-k}-1\\
\end{vsmallmatrix} 
=& ((\beta/\alpha)^{l-k}-1)((\gamma/\alpha)^{m-k}-1) \\ 
\nonumber &- ((\beta/\alpha)^{m-k}-1)((\gamma/\alpha)^{l-k}-1).
    \end{align*}
So the ratios $\beta/\alpha$ and $\gamma/\alpha$ lie on a complex algebraic curve;
thus we have more options than just roots of unity.

Geometrically, the vanishing of
    \begin{align*} 
 \begin{vsmallmatrix}
1& 1 & 1\\
1&(\beta/\alpha)^{l-k}  & (\beta/\alpha)^{m-k} \\
1&(\gamma/\alpha)^{l-k} & (\gamma/\alpha)^{m-k}\\
\end{vsmallmatrix} 
    \end{align*}
means that $(1,1)$, $((\beta/\alpha)^{l-k}, (\beta/\alpha)^{m-k})$ and $((\gamma/\alpha)^{l-k}, (\gamma/\alpha)^{m-k})$ are collinear in $\CC^2$.
Therefore, we can construct these points using secants of the monomial curve $t \mapsto (t^{l-k}, t^{m-k})$ that pass through $(1,1)$.

For generic choices of $\beta/\alpha$ and $\gamma/\alpha$, the rank-2 flats of $M$ are the sets of cardinality 2 and the shifts of the set $\aset{k,l,m}$.

Depending on the arithmetic properties of $l-k$ we might get larger flats, including infinite ones. 
We explore this possibility in the remaining example.

\subsubsection{Some periodicity without roots of unity}
If $\aset{k,l,m} = \aset{0, n, 2n}$, the $n$-th shift yields $\aset{n, 2n, 3n}$ as a dependency, so $\aset{0,n, 2n, 3n}$ is a dependency, and so on.
We obtain $n\NN$ and its shifts as flats.
Thus, we have $n$ rank-2 flats with infinite cardinality $n\NN$,  $n\NN + 1$, etc. 
We also have that any pair $\aset{i,j}$ such that $i \not \equiv j \pmod{n}$ is a rank-2 flat.

\subsection{$p$-adic valuation}
    \label{sub:padicvaluationOrderThree}
\begin{figure}[htbp]
    \centering
    \begin{tikzpicture}[scale=0.8, >=stealth]
        
        \draw[thick] (3, 5) -- (3, 7.5);
        \draw[thick] (3, 5) -- (6, 5);
        \draw[thick] (3, 5) -- (0.5, 2.5);

        \filldraw (3, 5) circle (2pt) node[below right] {$(u_3, u_2 + u_3)$};

        \draw[thick, fill=white] (3, 6) circle (2pt) node[left=3pt] {$(u_3, 2u_3)$};
        
    \end{tikzpicture}
    \caption{A tropical line in the $(b_1,b_2)$-plane. 
    Any initial condition in the complement extends uniquely to a solution in $\bigcap_{n \in \NN}\SolSet(\trop(d^nF))$. 
The line itself and the empty dot are the projections to the $(b_1,b_2)$-plane of the tropicalization of $E_2$ and of $E_3$, respectively.}
    \label{fig:Order3Stratification}
\end{figure}
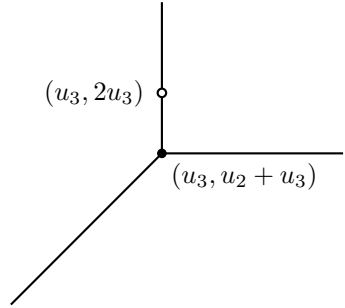

We extend the strategy used in the order-2 case, 
providing a conceptual explanation of the stable and unstable solutions.
The upshot is that if the initial conditions $(b_1,b_2)$ avoid a specific tropical line in the $(b_1,b_2)$-plane, then they determine a unique solution~$B$.
We will treat the realizability question of $B$ in an upcoming expansion of this note, 
where we also generalize our observations to arbitrary order.
We consider the differential polynomial
\begin{align*} 
    F &= (x-\alpha_1)(x- \alpha_2)(x - \alpha_3) \\ &= 
 x^3 - (\alpha_1 + \alpha_2 + \alpha_3)x^2 + (\alpha_1 \alpha_2 + \alpha_2 \alpha_3 + \alpha_1 \alpha_3)x - \alpha_1 \alpha_2 \alpha_3 \in \CC_p[\![t]\!]\{x\}. 
\end{align*}
The classical solutions are $\exp(\alpha_1t)$, $\exp(\alpha_2t)$ and $\exp(\alpha_3t)$.
We set $u_i = \nu_p(\alpha_i)$ for the valuations. As before we assume they are strictly ordered $u_1 < u_2 < u_3$.
We filter the space of solutions by their radii of convergence:
\begin{align*}
    E_1 &= \operatorname{span}( \exp(\alpha_1 t), \exp(\alpha_2 t), \exp(\alpha_3 t) ), \\
    E_2 &= \operatorname{span}( \exp(\alpha_2 t), \exp(\alpha_3 t) ), \\
    E_3 &= \operatorname{span}( \exp(\alpha_3 t) ).
\end{align*}
Indeed, $E_3 \subset E_2 \subset E_1$ 
and the generic element of $E_i$ has a radius of convergence equal to~$\rho(\exp(\alpha_i t))$.
Tropicalizing and then projecting these three spaces yields a stratification of the space of tropical initial conditions $b_0$, $b_1$, and $b_2$, namely $\trop(E_3) \subset \trop(E_2) \subset \trop(E_1)$.
See Figure~\ref{fig:Order3Stratification}, where we assume that $b_0 = 0$ and plot the $(b_1,b_2)$-plane.
To describe equations for this stratification, 
note that $E_1$, $E_2$ and $E_3$ are solution spaces to the differential polynomials:
\begin{align*}
    F_1 &= F, \\
    F_2 &= (x-\alpha_2)(x-\alpha_3),  \\
    F_3 &=  x-\alpha_3,
\end{align*}
respectively. Now we tropicalize, using the $u_1 < u_2 < u_3$ condition to determine the coefficients of the equations.
To streamline the exposition and highlight the combinatorics, 
we temporarily omit the factorial terms in the following equations. 
This corresponds to assuming that $p$ is sufficiently large relative to $i$:
\begin{align*}
        G^{(1)}_i &= \min(u_1 + u_2 + u_3 + b_i, u_1 + u_2 + b_{i+1}, u_1 + b_{i+2}, b_{i+3}), \\
    G^{(2)}_i &= \min(u_2 + u_3 + b_i, u_2 + b_{i+1}, b_{i+2})  \\
    G^{(3)}_i &=  \min(u_3 + b_i, b_{i+1}).
\end{align*}
Keeping track of factorial terms results in cancellations, as in the order 2 case. 

Given initial conditions $(b_1,b_2)$, we wish to compute $b_3$ from $G^{(1)}_0$.
Notice that the minimum of the first three terms in $G^{(1)}_0$ equals $u_1 + G^{(2)}_0 $.
The value of $b_3$ is determined if this minimum is achieved only once, 
namely $(b_1,b_2) \notin V^{\text{trop}}( G^{(2)}_0 )$.
In that case, $\boxed{b_3 = u_1 + G^{(2)}_0(b_1, b_2)}$.

Continuing with $(b_1,b_2) \notin V^{\text{trop}}( G^{(2)}_0 )$, 
for $b_4$ we consider 
\[ 
    G^{(1)}_1 = \min(u_1 + u_2 + u_3 + b_1, u_1 + u_2 + b_{2}, u_1 + b_{3}, b_{4}).
\]
We claim that among the first three terms, $u_1 + b_3$ is a unique minimum. Indeed:
\begin{itemize} 
    \item $u_1 + b_3 = 2 u_1 + G^{(2)}_0(b_1, b_2) \le 2 u_1 + u_2 + b_1 < u_1 + u_2 + u_3 + b_1 $
\item $u_1 + b_3 = 2 u_1 + G^{(2)}_0(b_1, b_2) \le 2 u_1 + b_2 < u_1 + u_2 + b_2 $
\end{itemize}
Thus, $b_4 = u_1 + b_3 = 2 u_1 + G^{(2)}_0(b_1,b_2)$. 
Continuing inductively, it can be proven that if $(b_1,b_2) \notin V^{\text{trop}}( G^{(2)}_0 )$, 
then in the calculation of $b_{n+1}$ the term containing $b_n$ is strictly smaller than all the others.
This logic, combined with keeping track of the factorial terms, yields the following theorem: 
\begin{thm} 
    \label{thm:Determined}
   Let $p$ be a prime and $ F = (x-\alpha_1)(x- \alpha_2)(x - \alpha_3)  \in \CC_p[\![t]\!]\{x\} $ a differential polynomial.
Set $u_i = \nu_p(\alpha_i)$ and $G^{(2)}_0 = \min(u_2 + u_3, u_2 + b_{1}, b_{2})$.  
Assume that $u_1  < u_2  < u_3 $. 
If the initial condition $(b_1,b_2)$ is not in $V^{\text{trop}}( G^{(2)}_0 )$, 
then the unique extension to a solution in $\bigcap_{n \in \NN}\SolSet(\trop(d^nF))$ is
\begin{align*} 
  b_n = (n-2) u_1 + G^{(2)}_0(b_1,b_2) - \nu_p(n!).
\end{align*}
\end{thm}

Regarding the other strata, a similar argument shows that if  $(b_1,b_2)$ is in $V^{\text{trop}}( G^{(2)}_0 )$ but not in $V^{\text{trop}}( G^{(3)}_0 )$, 
then there is a unique extension with radius of convergence $\rho(\exp(\alpha_2t))$.
If we allow the radius of convergence to change to $\rho(\exp(\alpha_1t))$, 
then a branching phenomenon happens again, as in the order 2 case.
Moreover, the solution with radius of convergence $\rho(\exp(\alpha_3t))$ is unique and equals $\trop(\exp(\alpha_3t))$. 
Its truncation to the $(b_1,b_2)$-plane is marked as an open point in Figure~\ref{fig:Order3Stratification}.

\section*{Acknowledgements}
The authors are grateful to the MPI MIS Leipzig, the University of Bologna, and the University of Warwick for their hospitality during several stages of this research.
S.M.~was partially supported by the Wallenberg AI, Autonomous Systems and Software Program (WASP) funded by the Knut and Alice Wallenberg Foundation.
A.V.~has received support from the Swiss National Science Foundation (grant no.\ 200142) and the UK Engineering and Physical Sciences Research Council (grant no.\ EP/X02752X/1).

\newcommand{\etalchar}[1]{$^{#1}$}
	
\vspace{0.5cm}

\end{document}